\documentclass[letterpaper, 10 pt, conference]{ieeeconf}  

\IEEEoverridecommandlockouts                              

\usepackage{graphics} 
\usepackage{epsfig} 
\usepackage{mathptmx} 
\usepackage{times} 
\usepackage{amsmath} 
\usepackage{amssymb}  
\usepackage{blindtext}
\usepackage{pgfplots}
\usepackage{mathtools, nccmath}
\pgfplotsset{compat=newest}
\usepackage{changepage}

\usepackage{mathtools}
\usepackage{color}
\usepackage{tikz}
\usetikzlibrary{angles, quotes}
\usetikzlibrary{patterns,hobby}
\usetikzlibrary{patterns.meta}
\usetikzlibrary{automata,positioning}
\usepackage[font=small,skip=0pt]{caption}
\usepackage{tikz-cd}
\usepackage{mathrsfs}
\let\labelindent\relax
\usepackage{enumitem}
\title{\LARGE \bf
How do we walk? Using hybrid holonomy to approximate non-holonomic systems*
}

\author{Maria Oprea$^{1}$ and William Clark$^{2}$
\thanks{*This work was funded by NSF grant DMS-1645643.}
\thanks{$^{1}$ M. Oprea is with the Center for Applied Mathematics, Cornell University, Ithaca, NY 14850, USA {\tt\small mao237@cornell.edu}}
\thanks{$^{2}$ W. Clark is with the Department of Mathematics, Cornell University, Ithaca, NY 14850, USA {\tt\small wac76@cornell.edu}}%
}

\newtheorem{defn}{Definition}
\newtheorem{prop}{Proposition}
\newtheorem{remark}{Remark}

\begin{document}

\setlength{\abovedisplayskip}{6pt}
\setlength{\belowdisplayskip}{6pt}

\maketitle
\thispagestyle{empty}
\pagestyle{empty}

\begin{abstract}

Why do we move forward when we walk? Our legs undergo periodic motion and thus possess no net change in position; however, our bodies do possess a net change in position and we are propelled forward. From a geometric perspective, this phenomenon of periodic input producing non-periodic output is holonomy. To obtain non-zero holonomy and propel forward, we must alternate which leg is in contact with the ground; a non-zero net motion can be obtained by concatenating arcs that would individually produce no net motion.
We develop a framework for computing the holonomy group of hybrid systems and analyze their behavior in the limit as the number of impacts goes to infinity.

\end{abstract}

\section{INTRODUCTION}

    Strategies for locomotion in nature are extremely diverse, ranging from walking to crawling, flying and swimming. However, at the core, the net displacement is always the result of the interaction between different internal body parts and the environment.
    
    A great example that illustrates this concept is the falling cat problem \cite{cat2}. If dropped upside down with zero angular momentum, the cat is able to right itself in the air and land on its feet by performing a series of internal shape changes. Although apparently violated, the angular momentum conservation law still holds if the rigid body assumption is dropped. A cat is indeed a complex system made out of different body elements, each having its own angular momentum. Using the model of two cylinders, joined in the middle, Kane and Scher \cite{cat2} managed to explain this phenomena, and showed that the nontrivial interaction between two identical parts of the body and the environment can lead to a global change in orientation. Later on, Montgomery  \cite{cat} managed to explain the phenomena elegantly, using the machinery of geometric phases, an approach we will be taking throughout this paper.

    Recent efforts have expanded the analysis to the more general case of undulatory locomotion \cite{undulating}, in which the environment acts on the system by imposing non-holonomic constraints. This is the situation for the motion of crawling terrestrial animals \cite{crawling}, some species of fish, as well as a large variety of artificial systems such as the snakeboard, the double wheeled robot \cite{undulating}, or the snake-like robot \cite{snake}. However, as noted in \cite{ruina}, the cases of realizable non-holonomic constraints are restricted to the situations of rolling without slipping, and skates, which are found in a limited number of locomotion strategies. We wish to extend the analysis of \cite{undulating} to situations in which different body parts can move independently to achieve a global displacement; the motion is allowed to be discontinuous and the constraints can be piece-wise holonomic. Although some progress in this direction has been achieved in \cite{piece_wise} for the leg-peg model, there is still a need for a unified theory of piece-wise holonomic systems.

  \begin{figure}
       \centering
       \includegraphics[width = 0.35\textwidth, angle = 270]{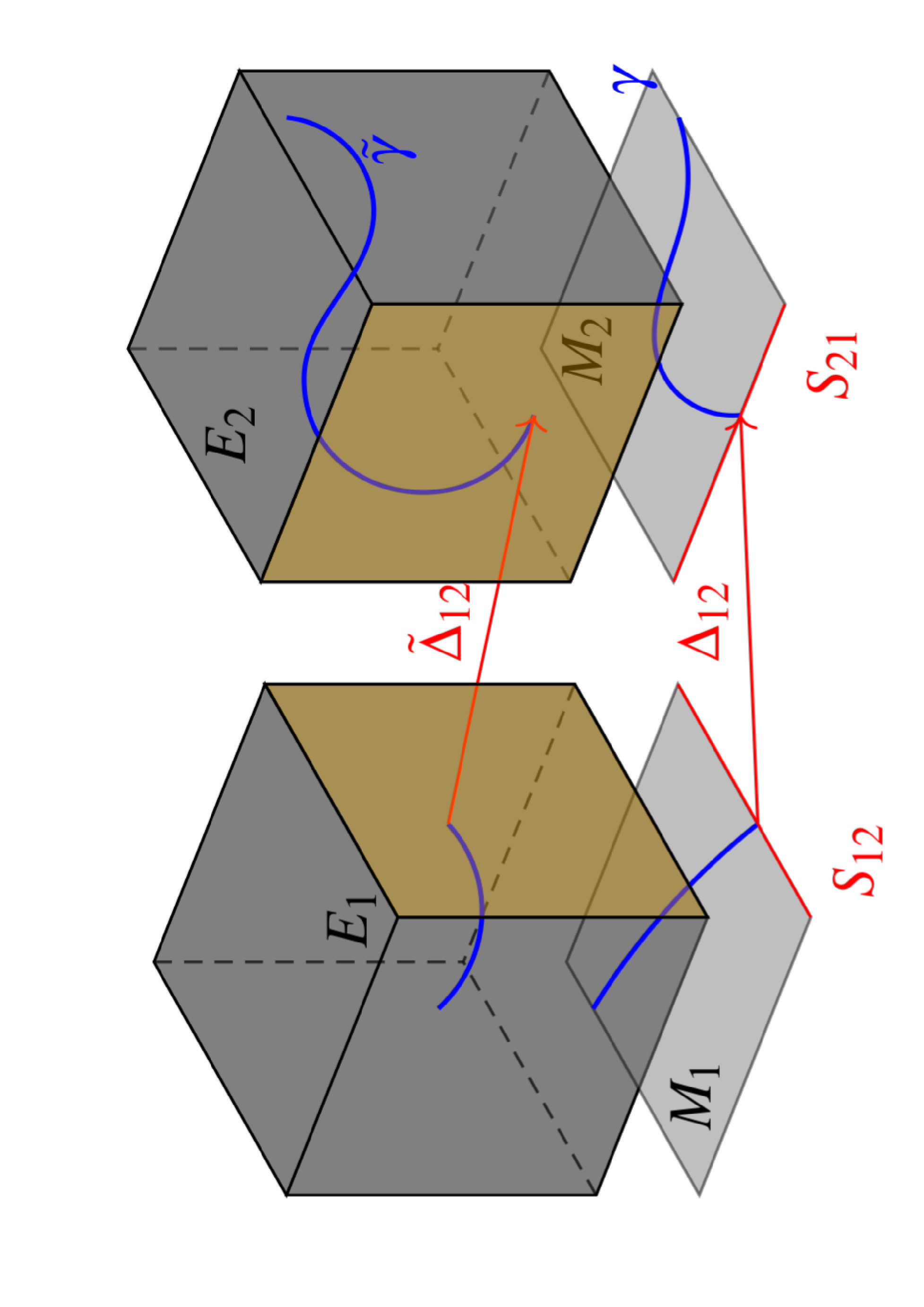}
       \setlength{\belowcaptionskip}{-15pt}
       \caption{Schematic picture of a hybrid bundle with two components. A curve $\gamma$ in the base space is lifted to $\tilde{\gamma}$ in the fiber space.}
      \label{fig:schematic}
   \end{figure}
  
    By considering piece-wise holonomic systems, we can accommodate for locomotion generated by walking, swimming or flying. In the case of walking: forward motion is obtained by switching the stance and swing legs, which undergo independent periodic motion. Although the body experiences a total net displacement, the feet do not "see" any change in their configuration between one body position and the other. Indeed the angle between the leg and the body, as well as its orientation relative to the hip is the same after performing step $n$ and step $n + 1$. Moreover, individually, each leg is subject to holonomic constraints, even though the entire body behaves like a non-holonomic object. In a similar fashion, the wings undergo the same periodic motion, their configuration does not change between flaps. Nevertheless, the bird is able to change altitude as well as horizontal position using just repeated changes in the state of its wings \cite{bird1}, \cite{bird2}. Likewise, the swimmer can generate forward movement using periodic motions of his arms and legs.

     However, in all the cases mentioned above, an isolated member is subject to holonomic constraints and cannot produce displacement on its own. Therefore, switching between different holonomic systems is essential in legged locomotion. Additionally, the cyclic motions of arms and legs are not smooth. The swimmer either takes his hand out of the water or bends it underneath the body in order to bring it back to the original position and the wings are straightened when pulled downwards and bent when the come back upwards. 
    In order to compensate for the holonomic nature of the constraint, the nature of the loops needs to be changed and switching between holonomic subsystems must be allowed. Consequently, we will obtain a system which is holonomic in each of its components, but non-holonomic as a whole. A framework for analyzing this interplay will be proposed in this paper.

    Even though motivation comes from locomotion, the model does not restrict to that. Examples such as a heaving buoy that uses an intermittent latching system to extract energy from the power of oscillating sea waves \cite{barazza},  or satellites that can reorient themselves using only internal rotors \cite{satellite} can also be interpreted in this framework.

   Since both piece-wise holonomic, and non-holonomic systems can generate net motion by manipulating their components, a natural question arises: is there a relationship between the two? In particular, can we approximate the smooth undulatory motion under non-holonomic constraints by a sequence of piece-wise holonomic systems, as we switch infinitely ofter between the components? The advantage of doing this is that we could transfer some of the machinery used for studying holonomic systems, which are Hamiltonian, to their nonholonomic counterparts. The idea of taking the limit in piece-wise holonomic systems was first presented by Andy Ruina \cite{ruina}, who conjectured that the motion of the smooth system and their non-smooth sister converge as some characteristic distance in the discrete system goes to zero. We formalise this and answer positively purely from a kinematic point of view without considering the dynamics of the system. In other words, the set of possible loops is not restricted to the integral curves of a vector field and our analysis is narrowed to the changes in space configuration and not in the velocities. 

  In our analysis we will employ the formalism of geometric phases \cite{geometric_phases}. As explained in \cite{intro}, the idea of internal shape changes generating total net movement can be mathematically expressed as a loop in a principle bundle generating a non-zero holonomy element in the fiber space. 
    
  More explicitly, we split he system variables into two parts: the one that describes the internal configuration  and the one representing the global status, called base space and fiber space respectively. Motion in the fiber space is equivalent to a net displacement, which is independent of the start and end point. This motivates the representation of net changes in the external configuration by elements of a Lie group, acting on the fiber space: $x_{final} = change(x_{initial}) = \Delta x (x_{initial})$.
    
  We consider a cyclic change in the internal configuration and ask ourselves: what is the effect of this motion in the external configuration? In order to connect the two we use the constraints and symmetries of the system, which give us a set of allowable directions of motion in the fiber space. Thus, the physics of the system will dictate a connection \cite{undulating,symmetry} which gives us a way of ``connecting'' motion in the base variables with that in the fiber variables. Our aim will be to compute the holonomy group element corresponding to a circular motion in the base space given by the net fiber displacement. However, unlike \cite{undulating} our base and fiber spaces will be composed out of discrete pieces, and multiple transitions are allowed to happen during a full cycle in the base. In order to deal with this we will make use of the hybrid system formalism.

In \S\ref{sec:prelim}, we offer an overview of the existing machinery of geometric mechanics, and extend it to the case of a hybrid system in \S\ref{sec:hybrid}. We then turn to the problem of computing the holonomy group of a hybrid principle bundle, and towards the end of \S\ref{sec:pb_formulation} we analyze the limit as the number of switches goes to infinity in the particular case of two alternating systems. Lastly, in \S\ref{sec:example}, we illustrate our theory by considering the double legged robot in one dimension.

\section{PRELIMINARIES}\label{sec:prelim}

We begin by defining the fundamental object of our formulation 
\begin{defn}[Principal bundle \cite{nonholo_book}, \cite{schuller}] Let $G$ be a Lie group and $E \xrightarrow{\pi} M$ a bundle. If:
\begin{enumerate}
\item $E$ is a right $G$ space i.e there exists a group action $ \bigstar:E \times G \rightarrow E$,
\item $\bigstar$ is free,
\item $E \xrightarrow{\pi} M$ and $E \xrightarrow{\rho} E/G$ are isomorphic as bundles where $\rho:E \rightarrow E/G$ is the natural projection onto the quotient space. 
\end{enumerate}
We call  $E \xrightarrow{\pi} M$ a principal $G$-bundle.
\end{defn}
Intuitively a principal bundle is a bundle whose fibers consist of Lie group elements. For us the total space $E$ will represent the configuration space of the system, while the base space $M$ will contain the variables we can control. Our overall goal will be to find an expression for the change in the fiber variables in terms of the base space variables. In this paper we will only consider $G$ abelian and isomorphic to $(\mathbb{R}^n, +)$.

\begin{remark}\label{rm:unique_g}
Going from $x$ to $y$ in the fibers corresponds to acting with a Lie group element. Moreover, since the action is free, that Lie group element will be unique i.e. for all $ x, y \in \pi^{-1}(m)$ there is a $g \in G $ such that $y = x \bigstar g$.
\end{remark}
We will show that no overall motion in the base space can generate nontrivial motion in the fiber space. In order to formalise this we consider a closed loop in the base space $\gamma:[0, 1] \rightarrow M$, with $\gamma(0) = \gamma(1) = m$. We are interested in what happens to the endpoints of the lifted loop $\tilde{\gamma}:[0, 1] \rightarrow E$ with $\pi(\tilde{\gamma}) = \gamma$. A priori, this is not well defined since there is no unique way of lifting a base curve to the fibers. We need a connection in order to have a systematic way of lifting $\gamma$. 

We call $ker(\pi_*)(m) = V_mE \subset T_mE$ the vertical space of the tangent bundle at $m$.
\begin{defn}(Principal connection)
A connection on a principal $G$ bundle is an assignment $m \mapsto H_mE$ such that:
\begin{enumerate}
\item $V_mE \oplus H_mE = T_mE$
\item $(\bigstar g)_*H_mE = H_{m \bigstar g}E$ (invariant  under the group action)
\item there is a unique decomposition $X_m = hor(X_m) \oplus vert(X_m)$ where $hor(X_m) \in H_mE$ and $vert(X_m) \in V_mE$.
\end{enumerate}
\end{defn}

Intuitively, a connection is a way to specify what a horizontal vector in the tangent space is. If we are given such an object we can define the horizontal lift of any curve $\gamma$ through some point $\tilde{\gamma}(0) \in \pi^{-1}(\gamma(0))$ as the unique curve $\tilde{\gamma}: [0, 1] \rightarrow E$ satisfying:
\begin{enumerate}
\item $\pi\circ\tilde{\gamma} = \gamma$
\item $\tilde{\gamma}'(x) \in H_{\gamma(x)}E $, for all $  x \in [0, 1]$ (all tangent vectors are horizontal)
\item $\pi_*(\tilde{\gamma}'(x)) = \gamma'(x)$,  for all $ x \in [0, 1]$ (the projected tangent vectors match the base tangent vectors).
\end{enumerate} 
Throughout this work we will assume that the connection is fixed and given by the physics of the problem, and we will try ask the question: what happens to the endpoint of the unique lifted loop? Note that since $\tilde{\gamma}(0)$, $\tilde{\gamma}(1) \in T_{\gamma(0)}E$ belong to the same fiber, by Remark \ref{rm:unique_g} there exists a $g$ such that $\tilde{\gamma}(1) = \tilde{\gamma}(0)\bigstar g$.
\begin{defn}[The holonomy group \cite{schuller}]
Let $m \in M$. Then the holonomy group at $m$ is 
\begin{multline}
Hol_m = \{g_{\gamma} | \tilde{\gamma}(1) = \tilde{\gamma}(0) \bigstar g_{\gamma} \text{ for some loop }\\\gamma:[0, 1] \rightarrow M, \ \gamma(0) = \gamma(1) = m\}.
\end{multline}
\end{defn}
The holonomy group is a subgroup of $G$ and it depends on the connection.

To compute the holonomy group it is useful to introduce the notion of a connection one form.  At a given point this form will give us the infinitesimal horizontal direction. 
\begin{defn}[Infinitesimal generators \cite{connection_curvature}]
Let $E \xrightarrow{\pi} M$ be a principal $G$ bundle. Then every element $A \in T_eG$, the Lie Algebra of $G$ induces a vertical vector field on $E$, $X^A$ defined as:
$$X^A f = f(m \bigstar exp(tA))'(0), \text{ for any } f \in \mathcal{C}^{\infty}(E). $$
The Lie Algebra element $A$ is called the infinitesimal generator of $X^A$.
\end{defn}
It can be shown that there is an isomorphism between the Lie Algebra of $G$ and the vertical vectors \cite{schuller}. Hence, for any vertical element $Y \in T_mE$ there is a unique infinitesimal generator $A$ such that $Y = X^A$.

\begin{defn}[Connection 1-form]
A connection 1-form, denoted by $\omega$, is a Lie algrebra valued form on the tangent bundle of the total space i.e. $\omega_m:T_mE \rightarrow T_eG, \ \text{for all} \ m \in M$ that returns for each vector in the tangent space the infinitesimal generator of its vertical part.
$$ \omega_m(Y) = A \ such \ that \ X^A = vert(Y), \ \text{for all} \ Y \in T_mE. $$
\end{defn}

In the following, we will provide a general procedure for computing the local holonomy group. Suppose we have a local section $\sigma: U \rightarrow \pi^{-1}(U)$. This induces a local trivialization $h: M\times G \rightarrow E$ which gives a local representation of the connection in the induced coordinates $h^*\omega$, as well as a projection of the connection in the base space variables $\sigma^*\omega: TM \rightarrow T_eG$. Suppose $\gamma: [0, 1]\rightarrow U$. Then the local holonomy group of an element $x \in \pi^{-1}\gamma(0)$ generated by the curve $\gamma$ is given by:
$$\Delta g_{loc} = \int_{\gamma}\sigma^* \omega.$$

We will denote the local coordinates by $(m, g)$ where $g = \sigma(m)$ Then the connection 1-form acting on tangent vectors $(\dot{m}, \dot{g})$ can be written as \cite{nonholo_book}:
\begin{equation}
\omega(\dot{m},\dot{g}) = (l_{g^{-1}})_* \dot{g} + A(m)\dot{m},
\end{equation}
where $l_{g^{-1}}$ denotes the left multiplication by the group element $g$ and $(l_{g^{-1}})_*$ is the lifted action on the tangent space. We need this in order to make sure that the result will be in $T_eG$ and not $T_gG$.

Thus we can represent $\omega$ in local coordinates as:
\begin{equation}\label{eq:local_repr}
 \omega = (Ad_{g})_*(l_{g^{-1}}^*dg + A(m)dm).
\end{equation}

Note that the adjoint is needed because the group $G$ acts on $E$ from the left and as such the allowed operator on $E$ is the right multiplication which is the adjoint of left multiplication.

\begin{prop}
Let $\gamma:[0, 1] \rightarrow U$ be a curve in the base space with $U \subset M$ the domain of the trivialization . Let $\tilde{\gamma}$ be its lift through $g_0 \in \pi^{-1} \gamma(0)$. Suppose we have a connection 1-form given in coordinates by \eqref{eq:local_repr}. Then the holonomy element can be computed as:
\begin{equation}\label{holonomy_computation}
    \Delta g_{loc} = \int_{\gamma}-A(m)dm.
\end{equation}

\end{prop}
\begin{remark}
 Expression \eqref{holonomy_computation} does not depend on the start and endpoint of the lift anymore. Indeed, the holonomy element only depends on the net change, and not the absolute position.
\end{remark}

We ilustrate the concepts by considering the rolling disk in 1D. The state space of the system is $(x, \theta)\in \mathbb{R}\times S^1$, where $x$ denotes the position fo the center of the ball on the the real line, and $\theta$ denotes the rotation angle at a given time. We can view this as a principal $\mathbb{R}-$bundle: $\mathbb{R}\times S^1 \xrightarrow{\pi} S^1 $ where $\mathbb{R}$ acts on $\mathbb{R}\times S^1 $ by translation in the fiber coordinate, i.e.
$ g\bigstar(x, \theta) = (x + g, \theta).$

Let $r$ be the radius of the ball. The motion constraints can be expressed as: $\dot{x} =r\dot{\theta}$. They naturally give rise to the connection: $Hor(X) = ker(\omega)$, where $\omega = dx - rd\theta$ is the connection 1-form. In order to compute the holonomy group we need to integrate $dx$ over the horizontal lift of closed curves. Let $\gamma:[0, 1]\rightarrow S^1$ be a closed loop on the circle with $\gamma(0) = \gamma(1) = \theta_0$. Let $\tilde{\gamma}: [0, 1] \rightarrow \mathbb{R}\times S^1$ be its horizontal lift through $\tilde{\gamma}(0) = x_0$. Then
$$ \Delta x = \int_{\tilde{\gamma}}dx = \int_0^1 \omega(\tilde{\gamma}'(t)) + rd\theta((\tilde{\gamma})'(t))dt.$$ The first term vanishes since $\tilde{\gamma}'(t)$ is horizontal for any $  t$ by construction. We are left with:
$$ \Delta x = r\int_0^1d\theta(\tilde{\gamma}'(t))dt = r\int_0^1 d\theta(\gamma') = r\int_{\gamma'}d\theta = 2\pi rn.$$ where $n$ is the number of times $\gamma$ wraps around the circle.

\section{HYBRID BUNDLES AND CONNECTIONS}\label{sec:hybrid}

In order to generalize the concept of holonomy to a hybrid system, we need to extend the notion of a principal bundle to the hybrid case. We begin by defining a hybrid bundle:
\begin{defn}(Hybrid bundle)
A hybrid bundle is a tuple $(E_i, \pi_i, M_i, S_{ij}, \Delta_{ij})$ where 
\begin{enumerate}
\item  $E_i\xrightarrow{\pi_i} M_i$ is a fiber bundle, where  $E_i$ is a smooth manifold called the state space, and $M_i$ is also a smooth manifold called the base.
\item  $S_{ij} \subset M_i$ is a submanifold of codimension 1, called the guard 
\item $D_{ij} = (\Delta_{ij}, \tilde{\Delta
}_{ij}) $ with $\Delta_{ij} : S_{ij} \rightarrow M_j$ and $\tilde{\Delta}_{ij}: \pi_i^{-1}(S_{ij}) \rightarrow \pi_j^{-1}(M_j) $ such that the following diagram commutes:
\begin{equation}
\begin{tikzcd}
             \pi_i^{-1}(S_{ij}) \arrow[r, "\tilde{\Delta}_{ij}"] \arrow[d, "\pi_i"] &  \pi_j^{-1}(M_j) \arrow[d, "\pi_j"] \\
           S_{ij} \arrow[r,"\Delta_{ij}"] & M_j
        \end{tikzcd}
\end{equation}
\end{enumerate} 
\end{defn}

A hybrid bundle is a hybrid principal bundle if each  $E_i\xrightarrow{\pi_i} M_i$ is a principal $G-$bundle for the Lie group $G$.

\begin{remark}
A loop in the base space of a hybrid bundle is a curve $\gamma:[0, 1] \rightarrow \cup M_i$ with $\gamma(0) = \gamma(1)$ satisfying: $\lim_{t \to t_0^+}\gamma(t) = \lim_{t \to t_0^-}\Delta_{ij}(\gamma(t))$, for any $ t_0 $ such that $\gamma(t_0) \in S_{ij}$ and $\gamma'$ is not tangent to $S_{ij}$ at $t_0$ (transversality) . In a hybrid system, a loop will be made out of separate pieces, each belonging to a different $M_i$. Similarly a loop in the state space is $\tilde{\gamma}:[0, 1]\rightarrow \cup E_i$ satisfying the regularity condition: $\lim_{t \to t_0^+}\tilde{\gamma}(t) = \lim_{t \to t_0^-}\tilde{\Delta}_{ij}(\tilde{\gamma}(t))$, for all $ t_0 $ such that $\gamma(t_0) \in \pi_i^{-1} (S_{ij})$. 

The transversality condition tells us that no impact happens when the curve is tangent to the impact surface. Suppose $S_{ij}$ is given by the level surface of some function $\eta_{ij}:M_i\rightarrow \mathbb{R}$. Then the conditions for a transition to happen at $t_0$ are: $\gamma(t_0) \in S_{ij}$ and $d\eta_{ij}(\gamma'(t_0)) \neq 0$. 

\end{remark}

There are two important particular cases to consider. 
\begin{enumerate}[label=(\roman*)]
\item Assume that all $M_i$ are embedded in a higher dimensional space $M$. Let $S_{ij} = S$ be a common subsurface. Then our bundle will become $(E_i, \pi_i, M_i, S, D_{ij})$ with $\Delta_{ij} = Id$ in the base variables. Loops in this type of systems will be continuous.

\item  Assume $E_i = E$, $M_i = M$ and $S_{ij} = S$, but $D_{ij} \neq Id$. In this case we loose continuity in the loops at the cost of the smoothness of the state space. 
\end{enumerate}

Let $\gamma:[0, 1]\rightarrow \cup M_i$ be a curve in the hybrid base space. Suppose there is a sequence $ (t_n)\in [0, 1]$ such that $\gamma(t_n) \in S_{ij}$ for some $i, j \in \{1, \dots, n\}$ and let $(i(n-1), i(n))$ be the indices of the spaces $\gamma$ transitions from and to when passing through $S_{ij}$.The finite sequence $(t_n)$ will partition the unit interval into pieces whose image lies within the $M_{i(n)}$ component of the hybrid bundle. 
\begin{defn}[Hybrid lift]
$\tilde{\gamma}:[0, 1] \rightarrow \cup E_i$ such that:
\begin{enumerate}
\item $\tilde{\gamma}(0) = e_0$
\item $\tilde{\gamma}(t) = \tilde{\gamma}_{i(n)}(t)$ if $t \in [t_{n}, t_{n + 1}]$, where $\tilde{\gamma_i}$ is the lift of $\gamma\cap M_i$ in the $E_i \xrightarrow{\pi_i} M_i$ bundle passing through $\tilde{\Delta}_{i(n-1)i(n)}(\tilde{\gamma}(t_n))$
\item $\tilde{\gamma}(t) = \lim_{t \to t_n^-} \tilde{\gamma}_{i(n - 1)}$ if $t \in S_{i(n - 1)}$\\
is called the hybrid lift of $\gamma$ passing through $e_0$.
\end{enumerate}
\end{defn}

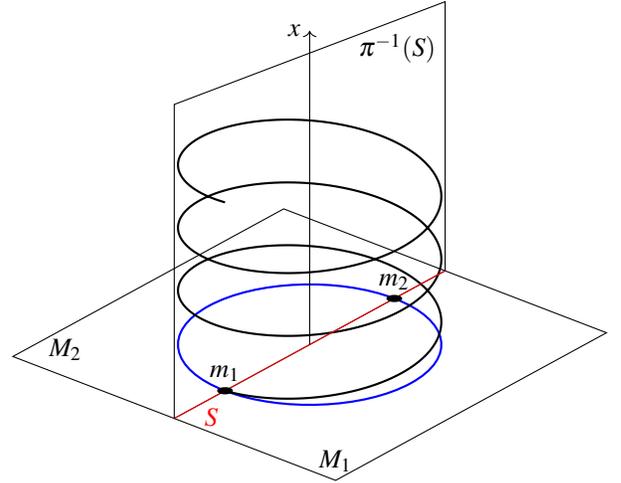
\begin{figure}
\centering
    \begin{tikzpicture}
    \begin{axis}[
 view={-40}{-40},
 axis line style = {draw=none},
 axis lines=middle,
 zmax=80,
  xmax=2,
   ymax=2,
 height=12cm,
 xtick=\empty,
 ytick=\empty,
 ztick=\empty,
 clip=false,
 x label style={draw = none},
   xlabel={$ $},
 y label style={at={(axis cs:0.05,2)},anchor=north},
   ylabel={$ $},
 z label style={at={(axis cs:0.075,0,80)},anchor=north},
   zlabel={$ $},
]
    \draw[thick, blue] (0,0) circle[radius=0.5];
    \addplot3+[domain=0:6*pi,samples=500,samples y=0,black,no marks, thick,
] 
({0.5*sin(deg(x))}, 
{0.5*cos(deg(x))}, 
{7*x/(pi)});
\draw (-0.8, -0.8, 0) -- (-0.8, 0.8, 0) -- (0.8, 0.8, 0) -- (0.8, -0.8, 0) -- cycle;
\draw (0, -0.8, 0) -- (0, 0.8, 0) -- (0, 0.8, 70) -- (0, -0.8, 60) -- cycle;
\node[above] at (0.8,0.8, 0) {$M_1$};
\node[below left] at (0,-0.8, 55) {$\pi^{-1}(S)$};

\node[label={[xshift=0.7cm, yshift=-0.3cm] $M_2$}] at (-0.8,0.8, 0){ };
\draw[red] (0,-0.8, 0) -- (0, 0.8, 0);
\node[red] at (0.1, 0.7, 0) {$S$};
\fill  (0,-0.5, 0) coordinate[label = above:$m_2$] circle(0.03);
\fill  (0,0.5, 0) coordinate[label = above:$m_1$] circle(0.03);

\draw [->] (0, 0, 0) -- (0, 0, 70) node[left]{$x$};
\end{axis}
    \end{tikzpicture}
    \setlength{\belowcaptionskip}{-15pt}
    \caption{Schematic drawing of a hybrid bundle ilustrating case (i) with two components $M_1$ and $M_2$ and the impact surface $S$. The  blue circle in the base space is lifted to the fiber space.}
\end{figure}

\begin{remark}
Suppose we have differential forms $\omega_i$ on each of the $E_i$ and a loop $\tilde{\gamma}$ in the state space with partition $(t_N)$ such that $t_0 = 0$ and $t_N= 1$. Then we define the integral of the tuple $(\omega_1, \dots, \omega_n)$ in the following way:
\begin{equation}
\int_{\tilde{\gamma}}(\omega_1, \dots, \omega_n) = \sum_{k = 0}^N \int_{t_k}^{t_{k + 1}}\omega_{i(k)}(\tilde{\gamma}'(t))dt.
\end{equation}
Note that each summand is a well defined integral performed in the space $E_{i(k)}$.
\end{remark}\label{integration}

\section{PROBLEM FORMULATION}\label{sec:pb_formulation}


Suppose we have a hybrid bundle $(E_i, \pi_i, M_i, S_{ij}, D_{ij})$ where each $E_i \xrightarrow{\pi_i} M_i$ is a principal $G$-bundle. Moreover, assume that each bundle $(E_i, \pi_i, M_i)$ is equipped with a connection 1-form $\omega_i: TE_i \rightarrow T_eG$ and a local section $\sigma_i:U_i \rightarrow E_i$ which induces a local trivialization $(m_i, g_i)$. Hence the connections in $U_i$ can be represented in coordinates as:
$$ \omega_i = Ad_{g_i}(l_{g_i^{-1}}^*dg_i + A(m_i, g_i)dm_i).$$
Further assume that $S_{ij} \cap U_i \neq \emptyset$ and that the form $Adm_i$ is exact in each of the $M_i$ i.e. there is a function $ F_i:M_i \rightarrow G$ such that $Adm = dF_i$.  If this is the case then the holonomy group af each of the $E_i \xrightarrow{\pi_i} M_i$ will be trivial which implies that no motion in the fibers is possible by cyclicly changing the variables in the base space. Now we want to study what happens to the holonomy group as we glue together these bundles. 

Consider a loop $\gamma:[0, 1]\rightarrow \cup U_i$ such that $\gamma(0) = \gamma(1)= m_0$. Let $(t_k)_{k = 1}^m$ be the times for which $\gamma(t_k) \in S_j$ for some $j \in \{1, \dots , n\}$ and let $i(k)$ be the index of the space $\gamma$ belongs to before making the transition at time $t_k$. Then we use Proposition \ref{holonomy_computation} and Remark \ref{integration} to compute the holonomy group:
\begin{align*}
\Delta g &= \int_{\gamma}(A_1 dm_1, \dots, A_n dm_n) = \sum_{k = 0}^N \int_{t_k}^{t_{k + 1}}A dm_{i(k)}(\gamma'(t))dt  \\
& = \sum_{k = 0}^m\int_{t_k}^{t_{k + 1}}dF_{i(k)}(\gamma'(t))dt  \\
& = \sum_{k = 0}^m F_{i(k)}(\gamma(t_{k + 1}^-)) - F_{i(k)}(\gamma(t_k)^+).
\end{align*}
We use the notation $\gamma(t_k)^+$ to write the value of the loop after the transition and $\gamma(t_k)^-$ to write the value before the  transition. They are related through the following identity:
$$\gamma^+(t_k) = \Delta_{i(k)i(k + 1)}\gamma^-(t_{k}).$$
This formula shows that the holonomy group only depends on what happens on the connection from which we obtain $F_i$ and what happens on the guard. We can write the holnomy group of a hybrid bundle as:
\begin{equation}\label{eq:final_holo}
\Delta g\left(\omega_i, \gamma \Big|_{\cup S_{ij}}\right) =\sum_{k = 0}^N F_{i(k)}(\gamma(t_{k + 1}^-)) - F_{i(k)}(\gamma(t_k)^+).
\end{equation}

We are now interested in the limit as the number of switches goes to infinity. In order to define the notion of a limit we assume that the $M_i$ are embedded into higher dimensional manifold and that $S_{ij} = S, \  i \in \{1, \dots, n\}, \ D_{ij} = Id$, as in case (i). Note however that the connection 1-form will change as we switch from one system to the other. Since the integrals do not depend of the curve parametrization assume that $(t_k)$ form an equidistant partition of $[0, 1]$ and call $\gamma(t_k) = m_k, \ k \in \{1, \dots, m\}$. In this paper we will consider the case $n = 2$,  the more general case will be presented in further work. Then equation \eqref{eq:final_holo} simplifies to:
\begin{equation}\label{eq:discrete_holo}
    \Delta g = \sum_{k = 0}^{N}\sum_{j = 1, 2}F_j(m_{k + 1}) - F_j(m_k).
\end{equation}

In particular let us consider a curve passing through points $m_1$ and $m_2$ in $S$ and making $N$ full loops. The summand will not depend on $k$ anymore, hence:
\begin{equation}\label{eq:2_step}
    \Delta g = N((F_1 - F_2)(m_2) - (F_1 - F_2)(m_1)),
\end{equation}
 where $(F_1 - F_2):S \rightarrow G$ continuous and differentiable.

Suppose there is an induced metric on $S$ and let $m_1 \to m_2$ in this metric. Then for a finite number of steps $\Delta g \mapsto 0$. However if we take infinitely many steps we obtain a nontrivial holonomy! 

Let $N \|m_1 - m_2 \|\rightarrow C$ for some constant $C$ as $N \to \infty$ and $m_1 \to m_2$. Then we get the infinitesimal version of the holonomy:

\begin{align*}
dg &= \lim_{m_1 \to m_2}\frac{(F_1 - F_2)(m_2) - (F_1 - F_2)(m_1)}{\|m1 - m2\|}\|m_1 -m_2\|N \\
&= Cd(F_1 - F_2)\Big|_{m_1}.
\end{align*}
\begin{remark}
 Our computations hold locally. However, we are interested in taking the limit as the number of impacts becomes infinite. This implies that the loops will get smaller and smaller, and the impact points will get closer together such that from some $m$ onwards, we can assume we are in the domain of a local section. 
 Since every closed form is locally exact, it is enough to assume that $Adm$ is closed, and restrict ourselves to $V :=U\cap U_{exact}$ where $U_{exact}$ is the neighbourhood on which $Adm $ can be written as $dF$ for some function $F: V\rightarrow G$.
\end{remark}

\section{EXAMPLE: PLANAR WALKER}\label{sec:example}

\begin{figure}
\centering
\begin{tikzpicture}
\draw[thick, ->] (-0.2,0) -- (5.2,0) coordinate[label=below: $\mathbb{R}$] (x);
\fill           (2.5,3) coordinate[label=above right:Hip] (z) circle(0.05);
\draw[dashed]  (2.5, 3)--(2.5,0) node[below] {$x$};
\coordinate[label=below left:$ O $] (O);
\fill           (5,0) coordinate (l) circle(0.05);
\draw[thick]    (O) to  (2.5,3);
\draw[thick]    (5, 0) to  (2.5,3);
\draw (O) -- (0, 1.5) coordinate (t);
\draw (5, 0) -- (5, 1.5) coordinate (p);
\pic [draw, <->,
      angle radius=11mm, angle eccentricity=1.2,
      "$\theta$"] {angle = z--O--t};
\pic [draw, <->,
      angle radius=11mm, angle eccentricity=1.2,
      "$\phi$"] {angle = p--l--z};
\pic [draw, <->,
      angle radius=11mm, angle eccentricity=1.2,left,
      "$2\delta$"] {angle = O--z--l};

\end{tikzpicture}
\setlength{\belowcaptionskip}{-15pt}
\caption{Diagram of the one dimensional walking robot in the moment when both feet are on the ground. The angle between the feet is the critical angle $2\delta$.}
\end{figure}
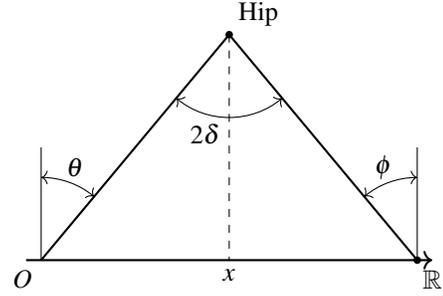

We will analyze the example of a planar walking robot with two legs. The state space variables are $x$, the horizontal position of the hip, and the angles $\theta$ and $\phi$ representing the angle of the two legs with respect to the vertical axis. As such, the state space is given by $\mathbb{R}\times S^1\times S^1$ which is a principal $\mathbb{R}$ bundle over $S^1 \times S^1$. The legs are two identical copies of each-other, which interact through the hip. Hence, the natural way to describe the walking robot is through a hybrid bundle $(E_i, \pi_i, M_i, S)$ $i\in\{1, 2\}$, where $E_1 = E_2 = \mathbb{R}\times S^1$ are the possible $x$ and stance angles,  $M_1 = S^1$ the possible $\theta$ and $M_2 = S^1$ the possible $\phi$ values. The impacts happen when both feet are on the ground, and the angle between the feet is $2\delta$, where $\delta$ is a characteristic of the robot. Thus, on the guard we have $\theta - \phi = 2\delta$ and $\theta = -\phi$. 

 When foot $\theta$ is down, the equations of motion $\dot{x} = l\cos{\theta}\dot{\theta}$ give rise to the connection: $\omega_{\theta} = dx - l\cos{\theta}d\theta$, where $l$ is the length of the leg. Similarly we obtain connection $\omega_{\phi} = dx -l\cos{\phi}d\phi$ for foot $\phi$. Hence in coordinates $A_{\theta} = l\cos{\theta}$ and $A_{\phi} = l\cos{\phi}$ which generate the exact forms: $A_{\theta}d\theta$ and $A_{\phi}d\phi$.
 \begin{remark}\label{rm:types_of_loops}
  Given the structure of the legs \cite{robots} a full step on one leg corresponds to $\theta$ and $\phi$ going from $-\delta$ to $\delta$ and hence a loop in the base space will start at $-\delta$, go all the way to $\delta$ and then come back to $\delta$. However, note that once $\theta$ reaches $-\delta$, the impact occurs. Hence in this case $m_1 = -\delta$ and $m_2 = -\delta$.
 \end{remark}

In order to move forward, the robot can only control $\theta$ and $\phi$. Suppose the robot only had the $\theta$ foot and let us compute the amount of forward movement it can generate after performing one full loop $\gamma$ in the base space.
$$ \Delta x = \int_{\gamma} dx = \int_{\gamma}l\cos{\theta}d\theta = l\sin{\gamma(1)} - l\sin{\gamma(0)} = 0.$$
There can be no forward movement using only one leg. Consider now a loop $\gamma$ in the base space of the hybrid bundle. Since $l\cos{\theta}d\theta$ and $l\cos{\phi}d\phi$ are exact, the integrals will be independent of the particular choice of $\gamma$ and will only depend on the endpoints. Let 
$$ \theta(t) = \begin{cases}\delta -4\delta t , \ \text{if } t < \frac{1}{2}\\
-\delta + 4\delta\left(t - \frac{1}{2}\right), \ \text{if } t \geq \frac{1}{2}\end{cases}$$
and $\gamma(t) = \begin{cases}\theta(t), & \text{if } t\leq \frac{1}{2}\\
-\theta(t)), &\text{if } t >\frac{1}{2}\end{cases}$. Note that $\gamma(t) \in M_1 = M_{\theta}$ for $t \leq \frac{1}{2}$ and $\gamma(t) \in M_2 = M_{\phi}$ for $t > \frac{1}{2}$. The impacts happen when the curve changes its slope, for $|\theta| = \delta = |\phi|$. Then using \eqref{eq:2_step}, we find that the local holonomy is given by:
$$ \Delta x = 4lN\sin\delta.$$
Given Remark \ref{rm:types_of_loops}, taking the limit $m_1 \to m_2$ is equivalent to letting
$\delta \to 0$. As the maximal angle between the legs decreases, we need to take more steps to achieve the same horizontal motion. We will assume that the two quantities are relate through: $\delta N\to \frac{C}{2}$,  as $\delta \to 0$ and $N \to \infty$. We obtain:
\begin{equation}\label{eq:walking_robot}
    \lim_{\delta \to 0}\Delta x = \lim_{\delta \to 0} 4lN\frac{\sin{\delta}}{\delta}\delta = lC.
\end{equation}
Suppose we increase our ``speed'' and we can take twice as many steps in the same amount of time i.e. $\tilde{N}  =2N$. Then $\delta \tilde{m} \to 2\frac{C}{2}$. We can think about $C$ as a variable that depends on how fast we can change the angle from $-\delta$ to $\delta$. It is natural therefore to consider infinitesimal amounts of changes in our speed of rotation, which are encoded in $dC$. Using equation \eqref{eq:walking_robot} we can relate the infinitesimal changes in $x$ to the infinitesimal changes in speed:
$$ dx = ldC.$$
If we rename $\theta := C$, this is exactly the form generating of the holonomy group of the rolling disk of radius $l$, whose connection 1-form is $\omega_{disk} = dx - ld\theta$. So in the limit as the number of steps goes to infinity and angle between the legs goes to zero, the walking robot behaves like the rolling disk!
\begin{figure}
\centering
\begin{tikzpicture}
\begin{axis}[grid=none,
    xlabel={$t$},
    xmax=5.3,ymax=2.6,
    xmin = -0.3, ymin = -2.6,
    xtick = \empty,
    ytick = \empty,
          axis lines=middle,
          restrict y to domain=-7:12,
          enlargelimits]
\addplot[red, domain =0:1]  {2.5 - 5*x} ;
\addplot[blue,dashed, domain=0:1]  {5*x - 2.5} ;

\addplot[red,dashed, domain =1:2]  {-2.5 + 5*(x - 1)} ;
\addplot[blue,domain=1:2]  {-5*(x - 1)  + 2.5} ;

\addplot[blue,dashed, domain = 2:3]  {-2.5 + 5*(x - 2)} ;
\addplot[red, domain=2:3]  {-5*(x - 2)  + 2.5} ;

\addplot[red, dashed, domain =3:4]  {-2.5 + 5*(x - 3)} ;
\addplot[blue, domain=3:4]  {-5*(x - 3)  + 2.5} ;

\addplot[blue,dashed, domain =4:5]  {-2.5 + 5*(x - 4)} node[below right] {$\phi(t)$};
\addplot[red, domain=4:5]  {-5*(x - 4)  + 2.5} node[above right]{$\theta(t)$};

\addplot[dashed, domain = 0:10] {2.5};
\addplot[dashed, domain = 0:10] {-2.5};

\node[left] at (0, 2.5) {$\delta$};
\node[left] at (0, -2.5) {$-\delta$};

\end{axis}
\end{tikzpicture}
\setlength{\belowcaptionskip}{-15pt}
\caption{ Two loops in $S^1$, represented in blue and red are combined to create a loop in the base space of the hybrid bundle $(\mathbb{R}\times S^1, \mathbb{R}\times S^1$).  The hard lines represent the components of the loop. The impacts happen when $\theta(t) = \delta$.}
\end{figure}
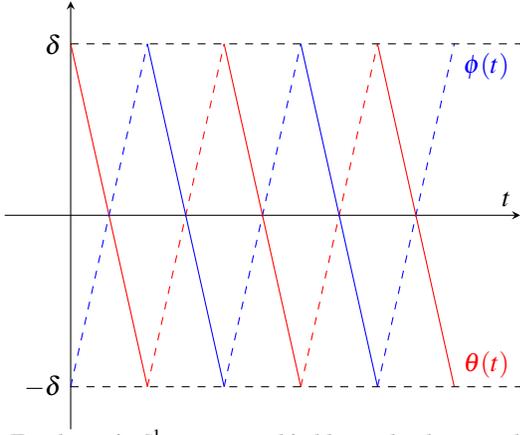
\section{CONCLUSIONS}

This paper builds up a framework for studying piece-wise holonomic systems by combining ideas from geometric mechanics as well as hybrid systems. Moreover, it offers a way of computing the local holonomy group of such a system, and analyzes what happens to the holonomy elements as the distance between two consecutive impacts goes to zero. 

There are multiple future research directions to consider. First of all, the work in this paper is restricted to an abelian group acting on the principal bundle, the behaviour of non-abelian groups is still to be analyzed. Moreover, the dynamics of the system were not considered, but in real life, the internal and external variables change according to some equations of motion that need to be accounted for in our formalism. Finally, we are interested in introducing controls to our system. In doing so we would be able to answer questions like: when is it optimal to transition from one piece to another in order to obtain the greatest amount of displacement in the shortest amount of time?

\addtolength{\textheight}{-12cm}   



\bibliographystyle{ieeetr}
\bibliography{references}

\end{document}